\documentclass[a4paper]{amsart}
\usepackage{amsmath,amsthm,amssymb}
\usepackage{graphicx}
\newtheorem{Def}{Definition}
\newtheorem{theorem}[Def]{Theorem}
\newtheorem{lemma}[Def]{Lemma}
\newtheorem{corollary}[Def]{Corollary}

\newtheorem*{problem}{Problem}
\makeatletter
 
\renewcommand{\@cite}[2]{{#1\if@tempswa , #2\fi}}
\makeatother

\title{The mapping class group of a punctured surface is generated by three elements}
\author{Naoyuki Monden}
\address{Department~of~Mathematics, Graduate~School~of~Science, Osaka~University, Toyonaka, Osaka~560-0043, Japan}
\email{n-monden@cr.math.sci.osaka-u.ac.jp}
\date{}
\begin{document}
\begin{abstract}
Let $\Sigma_{g,p}$ be a closed oriented surface of genus $g\geq 1$ with $p$ punctures. Let $\rm Mod(\Sigma_{\textit{g,p}})$ be the mapping class group of $\Sigma_{g,p}$.
Wajnryb proved in [Wa] that for $p=0, 1$ $\rm Mod({\Sigma_{\textit{g,p}}})$ is generated by two elements. Korkmaz proved in [Ko] that one of these generators can be taken as a Dehn twist.
For $p\geq 2$, We proved that $\rm Mod(\Sigma_{\textit{g,p}})$ is generated by three elements.  
\end{abstract}
\maketitle
\section{Introduction}
\ \ \ Let $\Sigma$ be a compact oriented surface of genus $g\geq 1$ with one boundary component, and $\Sigma_{g,p}$ be
a closed oriented surface of genus $g\geq 1$ with arbitrarily chosen $p$ points (which we call punctures). Let $\rm Mod(\Sigma)$ be the mapping class group of $\Sigma$, which is the group of homotopy classes of
orientation-preserving homeomorphisms which restrict to the identity on the boundary and let $\rm Mod(\Sigma_{\textit{g,p}})$ be the mapping class group of $\Sigma_{g,p}$, which is the group of homotopy classes of
orientation-preserving homeomorphisms which preserve the set of punctures. Let $\rm Mod^{\pm}(\Sigma_{\textit{g,p}})$ be the extended mapping class group of $\Sigma_{g,p}$, which is the group of homotopy class of all
(including orientation-reversing) homeomorphisms which preserve the set of punctures. By $\rm Mod_{\textit{g,p}}^{0}$ we will denote  the subgroup of $\rm Mod_{\textit{g,p}}$ which fixes the punctures pointwise. It is clear that we have
the exact sequence:
\begin{center}
$1 \rightarrow  \rm Mod_{\textit{g,p}}^{0} \rightarrow \rm Mod_{\textit{g,p}} \rightarrow Sym_{\textit{p}} \rightarrow 1,$
\end{center}
where the last projection is given by the restriction of a homeomorphism to its action on the puncture points.\\
\ \ \ The study of the generators for the mapping class group of a closed surface was first considered by Dehn. He proved in [De] that $\rm Mod(\Sigma_{\textit{g},0})$ is generated by a finite set of Dehn twists.
Thirty years later, Lickorish [Li] showed that $3\textit{g}-1$ Dehn twists generate $\rm Mod_{\textit{g},0}$. This number was improved to $2\textit{g}+1$ by Humphries [Hu]. Humphries proved, moreover, that in fact the number $2\textit{g}+1$ is minimal;
i.e. $\rm Mod(\Sigma_{\textit{g},0})$ cannot be generated by $2g$ (or less) Dehn twists. Johnson [Jo] proved that the $2g+1$ Dehn twists also generate $\rm Mod(\Sigma)$ and $\rm Mod(\Sigma_{\textit{g},1})$.
In the case of multiple punctures the mapping class group can be generated by $2g+p$ twists for $p\geq 1$ (see [Ge]).\\
\ \ \ It is possible to obtain smaller generating sets of $\rm Mod(\Sigma_{\textit{g,p}})$ by using elements other than twists. N.Lu (see [Lu]) constructed a generated  set of
$\rm Mod(\Sigma_{\textit{g},0})$ consisting of 3 elements. This result was improved by Wajnryb who found the smallest possible generating set of $\rm Mod(\Sigma_{\textit{g},0})$
consisting of 2 elements (see [Wa]). Korkmaz proved in [Ko] that one of these generators can be taken as a Dehn twist. It is also known that $\rm Mod(\Sigma_{\textit{g},0})$
can be generated by 3 torsion elements (see [BF]). More, Korkmaz showed in [Ko] that for $p=0, 1$ $\rm Mod(\Sigma_{\textit{g},{\textit{p}}})$ can be generated by 2 tosion elements.
\ \ \ In [Ma], Maclachlan proved that the moduli space is simply connected as a topological space by showing that $\rm Mod(\Sigma_{\textit{g},0})$ is generated by torsion elements.
Several years later Patterson generalized this result to $\rm Mod(\Sigma_{\textit{g,p}})$ for $g\geq 3$, $p\geq 1$ (see [Pa]).\\
\ \ \ In [MP], McCarthy and Papadopoulos proved that $\rm Mod(\Sigma_{\textit{g},0})$ is generated by infinitely
many conjugetes of a single involution for $g\geq 3$. Luo, see [Luo], described the finite set of involutions which generate $\rm Mod(\Sigma_{\textit{g,p}})$ for $g\geq 3$. 
He also proved that $\rm Mod(\Sigma_{\textit{g,p}})$ is generated by torsion elements in all cases except $g=2$ and $p=5k+4$, but this group is not generated by involutions if $g \leq 2$.
Brendle and Farb proved that $\rm Mod(\Sigma_{\textit{g,p}})$ can be generated by 6 involutions for $g\geq 3,p=0$ and $g\geq 4,p\leq 1$ (see [BF]). In [Ka], Kassabov proved that for every $p$
$\rm Mod(\Sigma_{\textit{g,p}})$ can be generated by 6 involutions if $g\geq 4$, 5 involutions if $g\geq 6$ and 4 involutions if $g\geq 8$. He also proved in the case of $\rm Mod^{\pm}(\Sigma_{\textit{g,p}})$.
In [St], Stukow proved that $\rm Mod(\Sigma_{\textit{g},0})$ is generated by 3 involutions for $g\geq 1$.\\
\ \ \ In this paper we show two results. We apply Korkmaz's result. We show that for $g\geq 1, p\geq 2$, $\rm Mod(\Sigma_{\textit{g,p}})$ is generated by 3 elements one of which is Dehn twists.
Next, we also prove that $\rm Mod^{\pm}(\Sigma_{\textit{g,p}})$ is generated by 3 elements one of which is Dehn twist.

\section{Preliminaries}
\ \ \ Let $c$ be simple closed curve on $\Sigma_{\textit{g,p}}$. Then the (right handed) Dehn twist $C$ about $c$ is the homotopy class of the homeomorphism obtained
by cutting $\Sigma_{\textit{g,p}}$ along $c$, twisting one of the side by $360^{\circ}$ to the right and gluing two sides of c back to each ohter. We denote
curves on $\Sigma_{\textit{g,p}}$ by letters $a,b,c,d$ and corresponding Dehn twists about them by capital letters $A,B,C,D$. \\
\ \ \ A small regular neighborhood of an arc $s_{ij}$ joining two puncuters $x_{i}$ and $x_{j}$ of $\Sigma_{g,p}$ is denoted by $N(s_{ij})$. The (right hand) half twist along $s_{ij}$ is denoted by $H_{ij}$. 
It is a self-homeomorphism $H_{ij}$ supported in $N(s_{ij}\cup x_{i}\cup x_{j})$ so that $H_{ij}$ leaves $s_{ij}$ invariant and interchanges $x_{i},x_{j}$ and $H_{ij}^{2}$ is the right handed Dehn twist along $\partial N(s\cup x_{i}\cup x_{j})$.\\
\ \ \ If $F$ and $G$ are two homeomorphisms, then the composition $FG$ means that $G$ is applied first.\\
\ \ \ We define the curves $a_{i}$, $b$, $\delta$ and the arc $s_{1p}$ on $\Sigma_{g,p}$ as shown in Figure 1.\\
\ \ \ We recall the following lemmas.
\begin{lemma}
Let $c$ be a simple closed curve on $\Sigma_{g,p}$, let $F$ be a self-homeomorphism of $\Sigma_{g,p}$ and let $F(c)=d$. Then $FCF^{-1}=D^{r}$, where $r=\pm 1$ depending on whether F is orientation-preserving or orientation-reversing.
\end{lemma}
\begin{lemma}
Let $c$ and $d$ be two simple closed curves on $\Sigma_{\textit{g,p}}$. If $c$ is disjoint from $d$, then $CD=DC$.
\end{lemma}

\begin{figure}[htbp]
 \begin{center}
  \includegraphics*[width=10cm]{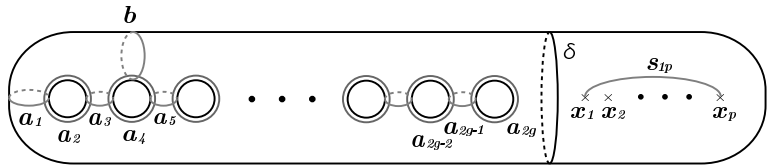}
 \end{center}
 \caption{The curves $a_{i}$, $b$, $c_{i}$, $d_{i}$}
 \label{fig:one}
\end{figure}
\newpage
\section{The mapping class group}
\ \ \ Let $S$ denote the product $A_{2g}A_{2g-1}\cdots A_{2}A_{1}$ of $2g$ Dehn twists in $\rm Mod(\Sigma_{\textit{g,p}})$ and let $G$ be the subgroup of $\rm Mod(\Sigma_{\textit{g,p}})$ generated by $B$ and $SH_{1p}$.
Lemma 3 follows the same argument as Section 3 of [Ko].
\begin{lemma}
$A_{1},\ldots ,A_{2g} \in G$.
\end{lemma}

Next, we prove that the subgroup $G^{\prime}$ generated by three elements of $\rm Mod(\Sigma_{\textit{g,p}})$ includes the pure mapping class group $\rm Mod^{0}(\Sigma_{\textit{g,p}})$, and that $G^{\prime}$ is equal to $\rm Mod(\Sigma_{\textit{g,p}})$.\\

\begin{figure}[htbp]
 \begin{center}
  \includegraphics[width=75mm]{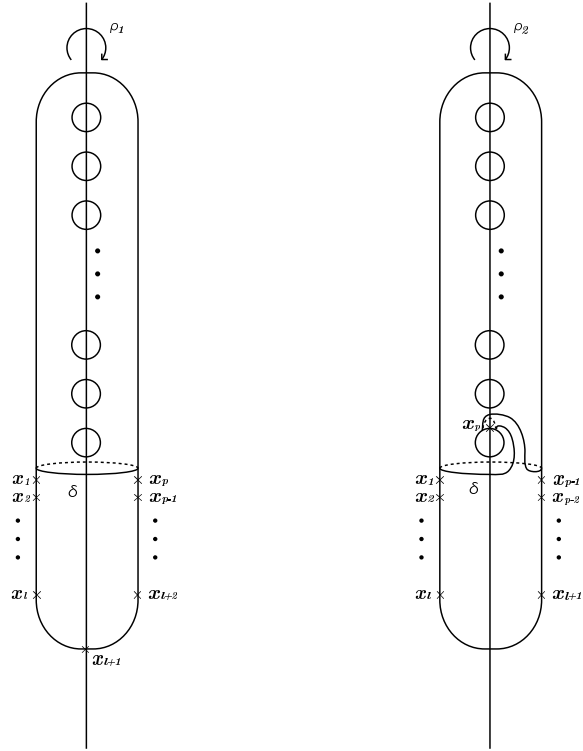}
 \end{center}
 \caption{The embeddings of the surface $\Sigma_{g,p}$ in the Euclidian space used to define the involutions $\rho_{1}$ and $\rho_{2}$.}
 \label{fig:one}
\end{figure}

Let us embed $\Sigma_{g,p}$ in Euclidian space in two different ways as shown on Figure 2. (In these pictures we will assume that the number of punctures $p$ is odd. In the case of even number of punctures we have to swap the lower parts of pictures.)
The curve $\delta$ separates $\Sigma_{g,p}$ into two components: the first one, denoted by $\Sigma$, is a surface of genus $g$ with boundary component and no punctures. The second one (denoted by $D$) is a disk with $b$ puncture points.\\
\ \ \ Each embedding gives a natural involution of the surface --- the half turn rotation around its axis of symmetry. Let us call these involutions $\rho_{1}$ and $\rho_{2}$, and denote the product $T=\rho_{1}\rho_{2}$.\\
\ \ \ On the set of punctures $T$ acts as a long cycle
\begin{center}
$T(x_{p})=x_{1}$ and $T(x_{i})=x_{i+1}$ for $1\leq i\leq p-1$.
\end{center}
\ \ \ Let $G^{\prime}$ be the subgroup of $\rm Mod(\Sigma_{\textit{g,p}})$ generated by $B$, $SH_{1p}$ and $T$. We prove that $G^{\prime}$ includes $\rm Mod^{0}(\Sigma_{\textit{g,p}})$.
In [Ge] it is shown that $\rm Mod^{0}(\Sigma_{\textit{g,p}})$ is generated by the Dehn twists around the curve $b$, $a_{i}$-es and $e_{j}$-es, for $j=1,\ldots ,p-1$, where the curves $e_{i}$-es are shown on Figure 3.

\begin{figure}[htbp]
 \begin{center}
  \includegraphics[width=30mm]{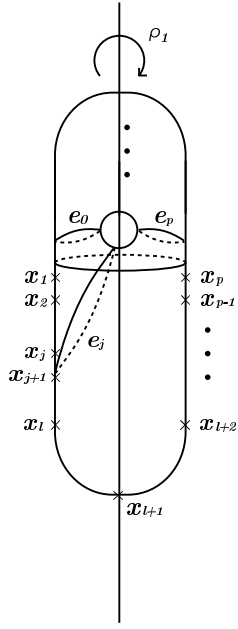}
 \end{center}
 \caption{The curves $e_{j}$.}
 \label{fig:one}
\end{figure}

\begin{lemma}
T acts on the curve $e_{j}$ as follows:
\begin{eqnarray*}
T(e_{i})=e_{i+1} \ \ (i=0,\ldots ,p-1).
\end{eqnarray*}
\end{lemma}
\begin{proof}
The Figure 4 shows the action of $\rho_{1}$ and $\rho_{2}$ on the curve $e_{i}$. It is clear from the picture that $e_{i+1}=\rho_{1}\rho_{2}(e_{i})=T(e_{i})$.
\end{proof}

\begin{figure}[htbp]
 \begin{center}
  \includegraphics[width=100mm]{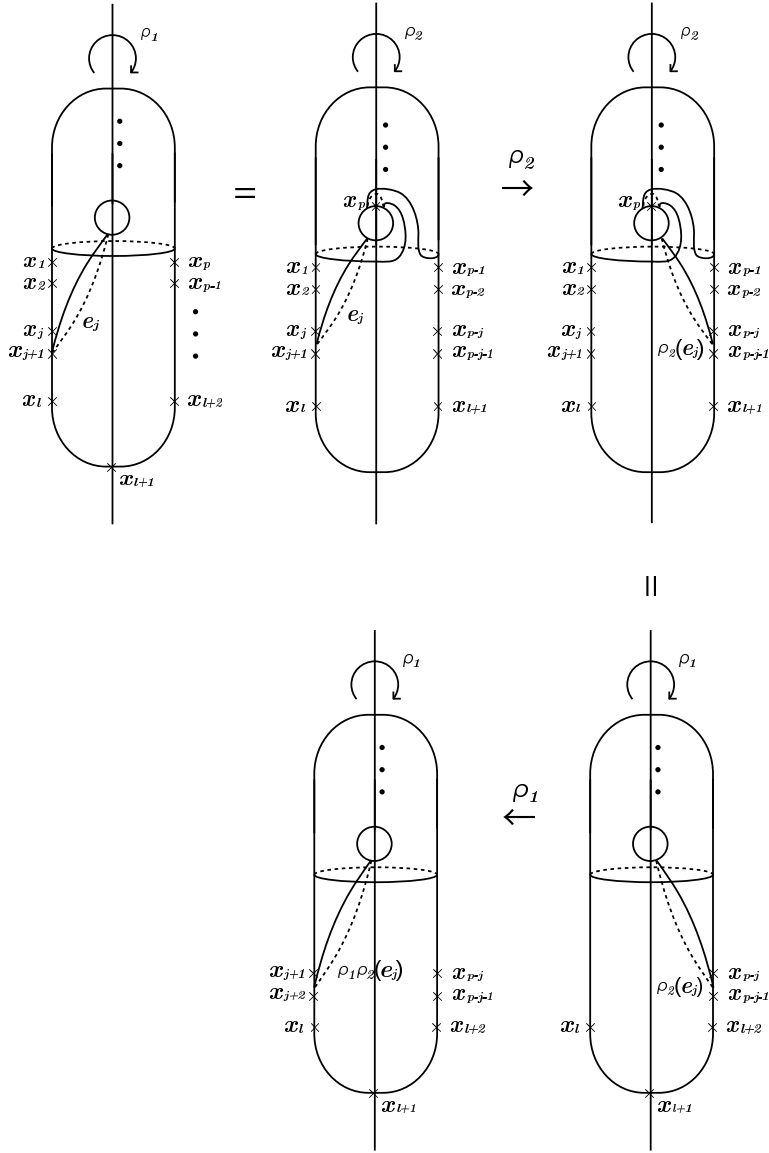}
 \end{center}
 \caption{The action of $T$ on the curve $e_{i}$.}
 \label{fig:one}
\end{figure}

\begin{lemma}
$E_{0},\ldots ,E_{p-1} \in G^{\prime}$
\end{lemma}
\begin{proof}
The inclusion $\iota :\Sigma \rightarrow \Sigma_{g,p}$ induces a natural mapping
\begin{eqnarray*}
\iota_{*}:\rm Mod(\Sigma)\rightarrow Mod(\Sigma_{\textit{g,p}}).
\end{eqnarray*}
If $F$ is in $\rm Mod(\Sigma)$, then $\iota(F)$ is represented by extending $F$ to $\Sigma_{g,p}$ using the identity mapping on $D$. In [Jo], Johnson proved that $\rm Mod(\Sigma)$ is generated by $B, A_{1},\ldots ,A_{2g}$.
Since $B, A_{1},\ldots ,A_{2g}$ are in $G\subset G^{\prime}$, $G^{\prime}$ contains $\iota_{*}(\rm Mod(\Sigma))$. Therefore $E_{0}$ is in $\iota_{*}(\rm Mod(\Sigma))\subset G^{\prime}$. Using Lemma 4 we can prove that all $E_{j}=T^{j}E_{0}T^{-j}$ are in $G^{\prime}$.
\end{proof}
\begin{corollary}
$\rm Mod^{0}(\Sigma_{\textit{g,p}}) \subset G^{\prime}$.
\end{corollary}
We prove that $G^{\prime}=\rm Mod(\Sigma_{\textit{g,p}})$. For it, we need the next lemma and corollary.
\begin{lemma}
Let $H$, $Q$ denote the groups and let $N, K$ denote the subgroups of $G$. Suppose that the group $H$ has the following exact sequence;
\begin{eqnarray*}
1 \rightarrow N \xrightarrow[]{i} H \xrightarrow[]{\pi} Q \rightarrow 1.
\end{eqnarray*}
If $K$ contains $i(N)$ and has a surjection to $Q$ then we have that $K=H$.
\end{lemma}
\begin{proof}
We suppose that there exists some $h\in H-K$. By the existence of surjection from $K$ to $Q$, we can see that there exists some $k\in K$ such that $\pi(k)=\pi(h)$.
Therefore, since $\pi(h^{-1}k)=\pi(h)^{-1}\pi(k)=1$, we can see that $h^{-1}k\in \rm Ker \ \pi= Im \ \textit{i}$. Then there exists some $n\in N$ such that $i(n)=h^{-1}k$.
By $i(N)\subset K$, since $i(n)\in K$ and $k\in K$, we have
\begin{eqnarray*}
h=k\cdot i(n)^{-1}\in K.
\end{eqnarray*}
This is contradiction in $h\notin K$. Therefore, we can prove that $K=H$.
\end{proof}
It is clear that we have the exact sequence:
\begin{eqnarray*}
1 \rightarrow  \rm Mod_{\textit{g,b}}^{0} \rightarrow \rm Mod_{\textit{g,b}} \rightarrow Sym_{\textit{b}} \rightarrow 1.
\end{eqnarray*}
Therefore, we can see the following corollary;
\begin{corollary}
Let $K$ denote the subgroup of $\rm Mod(\Sigma_{\textit{g,b}})$, which contains $\rm Mod^{0}(\Sigma_{\textit{g,b}})$ and has a surjection to $\rm Sym_{\textit{b}}$. Then K is equal to $\rm Mod(\Sigma_{\textit{g,b}})$.
\end{corollary}
\begin{theorem}
Suppose that $g\geq 1$ and $p\geq 2$. The subgroup $G^{\prime}$ generated by $B$, $SH_{1p}$ and $T$ is equal to the mapping class group $\rm Mod(\Sigma_{\textit{g,p}})$.
\end{theorem}
\begin{proof}
From Corollary 6, so that $G^{\prime}$ satisfies the condition of Corollary 8, we only need to show that $G^{\prime}$ can be made to map surjectively onto $\rm Sym_{\textit{p}}$.
This is equivalent to showing that the images to $\rm Sym_{\textit{p}}$ of the elements of $G^{\prime}$ generate $\rm Sym_{\textit{p}}$. It is clear that the image of $SH_{1p}$ is $(1,p)$ and the image of $T$
is $(1,\ldots ,p)$. Since $(1,p)$ and $(1,\ldots,p)$ generate $\rm Sym_{\textit{p}}$, we finish the proof of Theorem 9.
\end{proof}

\section{The extended mapping class group}
\ \ \ In this section we prove that the extended mapping class group $\rm Mod^{\pm }(\Sigma_{\textit{g,p}})$ is also generated by three elements.
In order to generate the extended mapping class group $\rm Mod^{\pm}(\Sigma_{\textit{g,p}})$, it suffices to add one more generator, namely the homotopy class of any orientation-reversing homeomorphism.\\

Let $R$ denote the reflection in Figure 7 and let $T^{\prime}$ denote the product $R\rho_{2}$.
Note that we replaced $\rho_{1}$ with $R$. Let $H$ denote the subgroup of $\rm Mod^{\pm }(\Sigma_{\textit{g,p}})$ generated by $B$, $SH_{1p}$ and $T^{\prime}$.

\begin{figure}[htbp]
 \begin{center}
  \includegraphics[width=70mm]{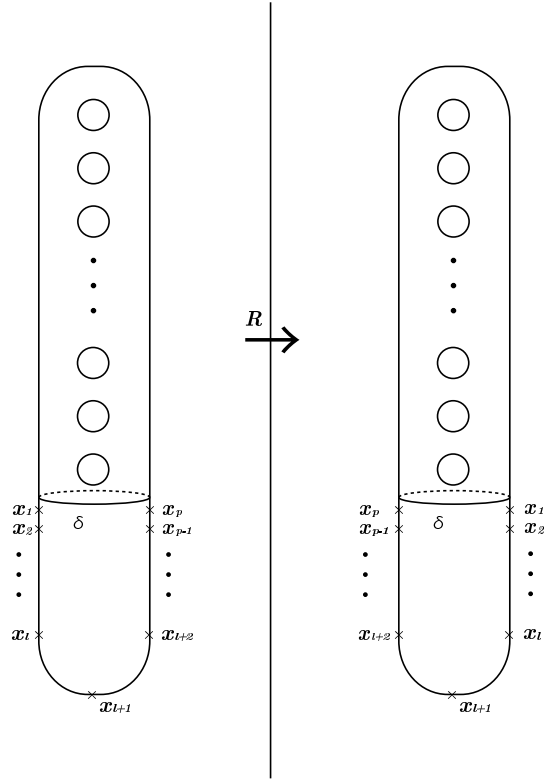}
 \end{center}
 \caption{The reflection $R$}
 \label{fig:one}
\end{figure}

\begin{theorem}
The subgroup $H$ generated by $B$, $SH_{1p}$ and $T^{\prime}$ is equal to the extended mapping class group $\rm Mod^{\pm }(\Sigma_{\textit{g,p}})$. 
\end{theorem}
\begin{proof}
By Section 3, 4, we can understand that $H$ contains $\rm Mod(\Sigma_{\textit{g,p}})$. Since $\rho_{2}$ and $T^{\prime}=R\rho_{2}$ are in $H$, $R$ is in $H$. We finish the proof of Theorem 10. 
\end{proof}

Wajnryb and Korkmaz proved in [Wa] and [Ko] that for $p=0, 1$ $\rm Mod(\Sigma_{\textit{g,p}})$ and $\rm Mod^{\pm }(\Sigma_{\textit{g,p}})$ can be generated by two elements. Then, we can consider following problem;

\begin{problem}
For any $p$ can $\rm Mod(\Sigma_{\textit{g,p}})$ and $\rm Mod^{\pm }(\Sigma_{\textit{g,p}})$ be generated by two elements?
\end{problem}

\end{document}